\documentstyle[11pt,A4]{article}
\def\mymod{{\hspace{10mm}\bmod\ }}
\def\QED{{\hfill$\Box$}}
\begin{document}
\bibliographystyle{plain}
\title{On the Periods of Generalized Fibonacci Recurrences$^1$}
\author{Richard P.\ Brent\\
	Computer Sciences Laboratory\\
	Australian National University\\
	Canberra, ACT~0200}
\date{TR-CS-92-03\\
	March 1992, revised March 1993}
\maketitle
\thispagestyle{empty}			%

\vspace{-1cm}
\begin{abstract}
We give a simple condition for a linear recurrence (mod $2^w$)
of degree $r$
to have the maximal possible period $2^{w-1}(2^r - 1)$.
It follows that the period is maximal in the cases of interest
for pseudo-random number generation,
i.e.\ for 3-term linear recurrences defined by trinomials
which are primitive (mod 2) and of degree $r > 2$.
We consider the enumeration of certain exceptional polynomials
which do not give maximal period, and list all such polynomials of
degree less than 15.

\end{abstract}

\section{Introduction}
\footnotetext[1]{{\em 1991 Mathematics Subject Classification.}
	Primary 11Y55, 		%
	12E05, 			%
	05A15; 			%
	Secondary 11-04, 	%
	11T06, 			%
	11T55, 			%
	12-04, 			%
	12E10, 			%
	65C10, 			%
	68R05.\\		%
\hspace*{16pt}{\em Key words and phrases.}
	Fibonacci sequence,
	generalized Fibonacci sequence,
	irreducible trinomial,
	linear recurrence,
	maximal period,
	periodic integer sequence,
	primitive trinomial,
	pseudo-random numbers.\\
\hspace*{16pt}Copyright \copyright\ 1992--2010, R.~P.~Brent.
	\hfill rpb133tr typeset using \LaTeX}
\addtocounter{footnote}{1}

The Fibonacci numbers satisfy a linear recurrence
\[F_n = F_{n-1} + F_{n-2} .\]
{\em Generalized Fibonacci} recurrences of the form
\begin{equation}
	x_n = \pm x_{n-s} \pm x_{n-r} \mymod 2^w 	\label{eq:GFR}
\end{equation}
are of interest because they are often used to generate
pseudo-random numbers~\cite{And90,Gre59,Jam90,Mar85a,Rei77,Tau65}.
We assume throughout that 		%
$x_0, \ldots, x_{r-1}$ are given and not all even,
and $w > 0$ is a fixed exponent.  Usually $w$ is close to the wordlength
of the (binary) computer used.

Apart from computational convenience, there is no reason to restrict
attention to 3-term recurrences of the special form~(\ref{eq:GFR}).
Thus, we consider a general linear recurrence
\begin{equation}
	q_0x_n + q_1x_{n+1} + \cdots + q_rx_{n+r} = 0 \mymod 2^w \label{eq:GLR}
\end{equation}
defined by a polynomial
\begin{equation}
	Q(t) = q_0 + q_1t + ... + q_rt^r 	\label{eq:Qdef}
\end{equation}
of degree $r > 0$.		%
We assume throughout that $q_0$ and $q_r$ are odd.
$q_0$ odd implies that the sequence $(x_n)$ is reversible,
i.e. $x_n$ is uniquely defined (mod $2^w$) by
$x_{n+1}, \ldots, x_{n+r}$.
Thus, $(x_n)$ is purely periodic~\cite{War33}.

In the following we often work in a ring $Z_m[t]/Q(t)$ of polynomials (mod $Q$)
whose coefficients are regarded as elements of $Z_m$ (the ring of integers
mod $m$).  For relations $A = B$ in $Z_m[t]/Q(t)$ we use the notation
\[	A = B \mymod (m, Q) .\]

It may be shown by induction on $n$ that if
$a_{n,0}, \ldots, a_{n,r-1}$ are defined by
\begin{equation}
	t^n = \sum_{j = 0}^{r-1} a_{n,j}t^j	\mymod (2^w, Q(t)) \label{eq:tn}
\end{equation}
then
\begin{equation}
	x_n = \sum_{j = 0}^{r-1} a_{n,j}x_j 	\mymod 2^w .	\label{eq:xn}
\end{equation}

Also, the generating function
\begin{equation}
	G(t) = \sum_{n \ge 0} x_nt^n 		\label{eq:G1}
\end{equation}
is given by
\begin{equation}
	G(t) = {P(t) \over \tilde Q(t)} \mymod 2^w, \label{eq:G2}
\end{equation}
where
\[ P(t) = \sum_{k=0}^{r-1} \left(\sum_{j=0}^{k} q_{r+j-k}x_j \right) t^k \]
is a polynomial
of degree less than $r$, and
\[	\tilde Q(t) = t^rQ(1/t) = q_0t^r + q_1t^{r-1} + ... + q_r \]
is the {\em reverse} of $Q$.
In the literature, $\tilde Q(t)$ is sometimes called the
{\em characteristic polynomial}~\cite{Gol67}
or the {\em associated polynomial}~\cite{War33}
of the sequence.
The use of generating functions is convenient and has been adopted by
many earlier authors (e.g.~Schur~\cite{Sch73}).	%
Ward~\cite{War33} does not explicitly use generating functions, but his
polynomial $U$ is the same as our $\tilde Q$, and many of his results
could be obtained via generating functions.

Let $\rho_w$ be the period of $t$
under multiplication mod $(2^w, Q(t))$,
i.e.\ $\rho_w$ is the least positive integer $\rho$ such that
\[ t^{\rho} = 1 \mymod (2^w, Q(t)) . \]	%
In the literature, $\rho_w$ is sometimes called the
{\em principal period}~\cite{War33} of the linear recurrence,
sometimes simply the {\em period}~\cite{Gol67}.	%
For brevity we define $\lambda = \rho_1$.

An irreducible polynomial in $Z_2[t]$
is a factor of $t^{2^r} - t$ (see e.g.~\cite{Wae49}),
so $\lambda \vert 2^r-1$.
We say that $Q(t)$ is {\em primitive} (mod 2) if $\lambda = 2^r - 1$.
Note that primitivity is a stronger condition than
irreducibility\footnote{For brevity we usually omit
the ``(mod 2)'' when saying that
a polynomial is irreducible or primitive.
Thus ``$Q(t)$ is irreducible (resp.\ primitive)'' means that $Q(t) \bmod 2$ is
irreducible (resp.\ primitive) in $Z_2[t]$.},
i.e.\ $Q(t)$ primitive implies that $Q(t)$ is irreducible,
but the converse is not generally true unless $2^r - 1$ is
prime\footnote{For example, the polynomial
$1 + t + t^2 + t^4 + t^6$ is irreducible, but not primitive, since it has
$\lambda = 21 < 2^6-1$.}. %
Tables of irreducible and primitive trinomials are
available~\cite{Gol67,Kur91,Rod68,Sta73,Wat62,Zie68,Zie69a,Zie69b,Zie70}.

In the following we usually assume that $Q(t)$ is irreducible.
Our assumption that $q_0$ and $q_r$ are odd excludes the trivial
case $Q(t) = t$, and implies that $\tilde Q(t)$ is irreducible (or primitive)
of degree $r$ iff the same is true of $Q(t)$.

We are interested in the period $p_w$ of the sequence $(x_n)$,
i.e.\ the minimal positive $p$ such that
\begin{equation}
	x_{n+p} = x_n \label{eq:period}
\end{equation}
for all sufficiently large $n$.
In fact, because of the reversibility of the sequence,
(\ref{eq:period}) should hold for all $n \ge 0$.
The period is sometimes called the {\em characteristic number}
of the sequence~\cite{War33}.
In general the period depends on the initial values
$x_0, \ldots, x_{r-1}$, but under our assumptions the
period depends only on $Q(t)$,
in fact $p_w = \rho_w$ (see Lemma~\ref{lemma:L_same_p}).

It is known~\cite{Knu81,Mar85b,War33} that
\[ p_w \le 2^{w-1}\lambda \]
with equality holding for all $w > 0$ iff it holds for $w = 3$.
The main aim of this paper is to give a simple necessary and sufficient
condition for
\begin{equation}
	p_w = 2^{w-1}\lambda . \label{eq:maxperiod}
\end{equation}
The result is stated in Theorem~\ref{thm:T1} in terms of a simple
condition which we call ``Condition~S''
(see Section~\ref{sec:cond_S}). %
In Theorem~\ref{thm:T2} we deduce that the period is maximal if
$Q(t)$ is a primitive trinomial of degree greater than 2.
Thus, in cases of practical interest for pseudo-random number
generation\footnote{A word of caution is appropriate.
Even when the period $p_w$ satisfies~(\ref{eq:maxperiod}),
it is not desirable to use a full cycle of $p_w$ numbers in applications
requiring independent pseudo-random numbers.  This is because only the
most significant bit has the full period. If the
bits are numbered from 1 (least significant) to
$w$ (most significant), then bit $k$ has period $p_k$.},
it is only necessary to verify that $Q(t)$ is primitive.
This is particularly easy if $2^r - 1$ is a Mersenne prime,
because then a necessary and sufficient condition is
\[	t^{2^r} = t	\mymod (2, Q(t)) .	\]

The basic results on linear recurrences modulo $m$ were
obtained many years ago~-- see for example Ward~\cite{War33}.
However, our main results (Theorems~\ref{thm:T1} and~\ref{thm:T2})
and the statement of ``Condition~S'' (Section~\ref{sec:cond_S})
appear to be new.

\section{A Condition for Maximal Period}
\label{sec:cond_S}

The following Lemma is a special case of
Hensel's Lemma~\cite{Knu81,Kri85,Zas69}
and may be proved using an application of Newton's method for
reciprocals~\cite{Kun74}.

\newtheorem{l_invert}{Lemma}	%
\begin{l_invert}
Suppose that $P(t) \bmod 2$ is invertible in $Z_2[t]/Q(t)$.
Then, for all $w \ge 1$, $P(t) \bmod 2^w$ is invertible in $Z_{2^w}[t]/Q(t)$.
\label{lemma:L_invert}
\end{l_invert}

We now give a sufficient condition for
the periods $p_w$ and $\rho_w$ to be the same.

\newtheorem{l_same_p}[l_invert]{Lemma}	%
\begin{l_same_p}
If $Q(t)$ is irreducible of degree $r$
and at least one of $x_0, \ldots, x_{r-1}$ is odd,
then $p_w = \rho_w$.
\label{lemma:L_same_p}
\end{l_same_p}

\leftline{\bf Proof}
\medskip

For brevity we write $p = p_w$ and $\rho = \rho_w$.
From (\ref{eq:G1}),
\[	G(t) = {R(t) \over 1 - t^p}	\mymod 2^w ,	\]
where $R(t)$ has degree less than $p$.
Thus, from (\ref{eq:G2}),
\begin{equation}
	R(t)\tilde Q(t) = (1 - t^p)P(t)	\mymod 2^w .	\label{eq:product1}
\end{equation}
Now $P(t) \bmod 2$ has degree less than $r$, but is not identically zero.
Since $\tilde Q(t) \bmod 2$ is irreducible of degree $r$,
application of the extended
Euclidean algorithm~\cite{Knu81} to $P(t) \bmod 2$ and $\tilde Q(t) \bmod 2$
constructs the inverse of $P(t) \bmod 2$ in $Z_2[t]/\tilde Q(t)$.
Thus, Lemma~\ref{lemma:L_invert} shows that
$P(t) \bmod 2^w$ is invertible in $Z_{2^w}[t]/\tilde Q(t)$.
It follows from~(\ref{eq:product1}) that
\[	t^p = 1			\mymod (2^w, \tilde Q(t)) ,\]
and $\rho \vert p$.
However, from (\ref{eq:tn}) and (\ref{eq:xn}), $p \vert \rho$.
Thus $p = \rho$.						 	\QED
\medskip

As an example, consider
$Q(t) = 1 - t + t^2$. %
We have $t^3 = 1 \bmod (2, Q(t))$,
$t^3 = -1 \bmod Q(t)$, and $t^6 = 1 \bmod Q(t)$, so
\begin{equation}
	\rho_w = \cases{3,	&if $w = 1$;\cr
	 		6,	&if $w > 1$.\cr}	\label{eq:rho3_6}
\end{equation}
It is easy to verify that~(\ref{eq:rho3_6}) gives the period $p_w$ of the
corresponding recurrence
\[	x_n = x_{n-1} - x_{n-2}	\mymod 2^w	\]
provided $x_0$ and $x_1$ are not both even.

The assumption of irreducibility in Lemma~\ref{lemma:L_same_p}
is significant.
For example\footnote{We thank a referee for suggesting this example.},
consider $Q(t) = t^2 - 1$ and $w = 1$, with initial
values $x_0 = x_1 = 1$. The recurrence is $x_n = x_{n-2} \bmod 2$,
so $p_1 = 1$, but $\rho_1 = 2$. Here $P(t) = 1 + t$ is a divisor of
$\tilde Q(t) = 1 - t^2$. %

We now define a condition which must be satisfied
by $Q(\pm t)$ if the period $p_w$
of the sequence $(x_n)$ is less than $2^{w-1}\lambda$
(see Theorem~\ref{thm:T1} for details).
For given $Q(t)$ the condition can be checked
in $O(r^2)$ operations\footnote{$O(r \log r)$ operations if the FFT is used
to compute the convolutions in~(\ref{eq:epssum}).}.
This is much faster than the method suggested by
Knuth~\cite{Knu81} or
Marsaglia and Tsay~\cite{Mar85b},	%
which involves forming high powers of $r \times r$ matrices (mod 8).

\medskip
\leftline{\bf Condition~S}
\medskip

Let $Q(t) = \sum_{j=0}^{r}q_jt^j$ be a polynomial of degree $r$.
We say that $Q(t)$ satisfies Condition~S if
\[	Q(t)^2 + Q(-t)^2 = 2q_rQ(t^2)	\mymod 8 .	\]

Lemma~\ref{lemma:L_cond_S} gives an equivalent condition\footnote{For another
equivalent condition, see~(\ref{eq:VW0}) and~(\ref{eq:VW5}).} which is more
convenient for computational purposes.  The proof is straightforward,
so is omitted.

\newtheorem{l_cond_S}[l_invert]{Lemma}
\begin{l_cond_S}
A polynomial $Q(t)$ of degree $r$ satisfies Condition~S iff
\begin{equation}
  \sum_{\scriptstyle j+k=2m \atop \scriptstyle 0 \le j < k \le r} q_jq_k
	= \epsilon_m	\mymod 2 		\label{eq:epssum}
\end{equation}
for $0 \le m \le r$, where
\begin{equation}
	\epsilon_m = {q_m(q_m - q_r) \over 2} .	\label{eq:epsdef}
\end{equation}
\label{lemma:L_cond_S}
\end{l_cond_S}

As an exercise, the reader may verify that the polynomial
$Q(t) = 1 - t + t^2$ satisfies both the definition of Condition~S
and the equivalent conditions of Lemma~\ref{lemma:L_cond_S}.
For other examples, see Table~\ref{Tab:exceptions}.

For convenience we collect some results regarding arithmetic in
the rings $Z_{2^w}[t]/Q(t)$.

\newtheorem{l_trivia}[l_invert]{Lemma}
\begin{l_trivia}
Let $X(t)$ and $Y(t)$ be polynomials over $Z$.
Then, for $w \ge 1$,
\begin{equation}
	X = Y \bmod (2^w, Q) \Rightarrow X^2 = Y^2 \bmod (2^{w+1}, Q) .
							\label{eq:trivia1}
\end{equation}
Also, if $Q(t)$ is irreducible, then
\begin{equation}
	X^2 = Y^2 \bmod (2, Q)	%
	\Leftrightarrow X^2 = Y^2 \bmod (4, Q)		\label{eq:trivia2}
\end{equation}
and
\begin{equation}
	X^2 = Y^2 \bmod (8, Q) \Leftrightarrow X = \pm Y \bmod (4, Q).
							\label{eq:trivia3}
\end{equation}
\label{lemma:L_trivia}
\end{l_trivia}

\leftline{\bf Proof}
\medskip
If $X = Y \bmod (2^w, Q)$ then $X = Y + 2^wR \bmod Q$ for some polynomial
$R(t)$ in $Z[t]$.  Thus

\noindent $X^2 = Y^2 + 2^{w+1}R(Y + 2^{w-1}R) \bmod Q$,
and~(\ref{eq:trivia1}) follows.

Now suppose that $Q(t)$ is irreducible.
If $X^2 = Y^2 \bmod (2, Q)$ then $(X - Y)^2 = 0 \bmod (2, Q)$.
Since $Q$ is irreducible, it follows that $X = Y \bmod (2, Q)$.
Thus, from~(\ref{eq:trivia1}), $X^2 = Y^2 \bmod (4, Q)$,
and~(\ref{eq:trivia2}) follows.

Finally, if $Q$ is irreducible and
$X^2 = Y^2 \bmod (8, Q)$ then,
as in the proof of~(\ref{eq:trivia2}), we obtain

\noindent $X = Y \bmod (2, Q)$, so
$X = Y + 2R \bmod Q$, where $R(t)$ is some polynomial in $Z[t]$.
Thus

\noindent $4R(Y + R) = 0 \bmod (8, Q)$,
i.e.~$R(Y+R) = 0 \bmod (2,Q)$.
Since $Q$ is irreducible, either

\noindent $R = 0 \bmod (2,Q)$ or
$Y + R = 0 \bmod (2,Q)$.  In the former case $X = Y \bmod (4,Q)$,
and in the latter case $X = -Y \bmod (4,Q)$.	%
Thus $X = \pm Y \bmod (4, Q)$. The implication in the other direction
follows from~(\ref{eq:trivia1}).	%
This establishes~(\ref{eq:trivia3}).
							\QED

The following Theorem is the key to the proof of Theorem~\ref{thm:T1}.
There is no obvious generalization to odd moduli.

\newtheorem{theorem0}{Theorem}
\begin{theorem0}
Let $Q(t) \bmod 2$ be irreducible in $Z_2[t]$.
Then
\[	t^\lambda = -1 \mymod (4, Q(t))	\]
iff $Q(t)$ satisfies Condition~S, and
\[	t^\lambda = 1 \mymod (4, Q(t))	\]
iff $Q(-t)$ satisfies Condition~S.
\label{thm:T0}
\end{theorem0}

\leftline{\bf Proof}
\medskip

Let
\[	V(t) = \sum_{j=0}^{\lfloor r/2 \rfloor} q_{2j}t^j,\;\;
	W(t) = \sum_{j=0}^{\lfloor(r-1)/2\rfloor} q_{2j+1}t^j ,\]
so $Q(t)$ splits into even and odd parts:
\begin{equation}
	Q(t) = V(t^2) + tW(t^2) .		\label{eq:VW0}
\end{equation}
By the definition of $\lambda$,
$t = t^{\lambda+1} \bmod (2, Q(t))$,
so
\begin{equation}
V(t^2) = t^{\lambda+1}W(t^2)	\mymod (2, Q(t)). \label{eq:VW1}
\end{equation}
Because $X(t^2) = X(t)^2 \bmod 2$ for any polynomial $X(t)$ in $Z[t]$,
(\ref{eq:VW1})~may be written as
\begin{equation}
	V(t)^2 = t^{\lambda+1}W(t)^2	\mymod (2, Q(t)).	\label{eq:VW05}
\end{equation}
$\lambda$, being a divisor of $2^r - 1$, is odd, so
$t^{\lambda+1}$ is a square.
Thus, from~(\ref{eq:trivia2}),
\begin{equation}
V(t)^2 = t^{\lambda+1}W(t)^2	\mymod (4, Q(t)). \label{eq:VW2}
\end{equation}
Also, since $V(t) = V(-t) \bmod 2$ and $W(t) = W(-t) \bmod 2$,
we have
\begin{equation}
V(-t)^2 = t^{\lambda+1}W(-t)^2	\mymod (4, Q(t)). \label{eq:VW25}
\end{equation}

To prove the first half of the Theorem, suppose that
\begin{equation}
	t^\lambda = -1 \mymod (4, Q(t)) .		\label{eq:minus_t}
\end{equation}
Thus, from~(\ref{eq:VW2}),
\begin{equation}
V(t)^2 + tW(t)^2 = 0	\mymod (4, Q(t)).		\label{eq:VW3}
\end{equation}
It follows that
\begin{equation}
	V(t)^2 + tW(t)^2 - q_rQ(t) = 0 \mymod (4, Q) .	\label{eq:VW4}
\end{equation}
However, the left hand side of~(\ref{eq:VW4}) is a polynomial of degree
less than $r$. Hence
\begin{equation}
	V(t)^2 + tW(t)^2 - q_rQ(t) = 0 \mymod 4 .	\label{eq:VW5}
\end{equation}
Replace $t$ by $t^2$ in the identity~(\ref{eq:VW5}).
From~(\ref{eq:VW0}), the result is easily seen to be equivalent to
$Q(t)$ satisfying Condition~S.

To prove the converse, suppose that $Q(t)$ satisfies Condition~S.
Reversing our argument,~(\ref{eq:VW3}) holds. Thus, from~(\ref{eq:VW2}),
\[	(t^{\lambda+1} + t)W(t)^2 = 0	\mymod (4, Q(t)).	\]
Now $W(t)$ has degree less than $r$, and $W(t) \ne 0 \bmod 2$
because otherwise, from~(\ref{eq:VW0}),
$Q(t) = V(t)^2 \bmod 2$ would contradict the irreducibility of $Q(t)$.
Thus, $W(t) \bmod 2$ is invertible in $Z_2[t]/Q(t)$.
From Lemma~\ref{lemma:L_invert}, $W(t) \bmod 4$ is invertible in
$Z_4[t]/Q(t)$, and we obtain
\[	t^{\lambda+1} + t = 0		\mymod (4, Q(t)).	\]
Since $Q(t) \ne t \bmod 2$, we can divide by $t$ to obtain
\[	t^\lambda = -1			\mymod (4, Q(t)).	\]
This completes the proof of the first half of the Theorem.

The proof of the second half is similar, with appropriate changes of sign.
Suppose that
\begin{equation}
	t^\lambda = 1 \mymod (4, Q(t)) .		\label{eq:plus_t}
\end{equation}
From~(\ref{eq:VW25}),
\begin{equation}
	V(-t)^2 = tW(-t)^2	\mymod (4, Q(t)) .	\label{eq:VW35}
\end{equation}
Thus, instead of~(\ref{eq:VW5}) we obtain
\begin{equation}
	V(-t)^2 - tW(-t)^2 - (-1)^{r}q_rQ(t) = 0 \mymod 4 . \label{eq:VW55}
\end{equation}
Replace $t$ by $-t^2$ in the identity~(\ref{eq:VW55}).
The result is equivalent to $Q(-t)$ satisfying Condition~S.
The converse also applies:
if $Q(-t)$ satisfies Condition~S then,
by reversing our argument and using irreducibility of $Q(t)$,
(\ref{eq:plus_t})~holds.						\QED
\medskip

We are now ready to state Theorem~\ref{thm:T1},
which relates the period of the sequence $(x_n)$ to Condition~S.
It is interesting to note that, in view of Theorem~\ref{thm:T0},
Theorem~\ref{thm:T1} is implicit in the discussion on page 628 of
Ward~\cite{War33}. More precisely,
Ward's case $T > 1$ corresponds to $Q(-t)$ satisfying Condition~S,
while Ward's case ($T = 1$, $K(x) = 1 \bmod 2$) corresponds to
$Q(t)$ satisfying Condition~S.
However, Ward's exposition is complicated by consideration of
odd prime power moduli (see for example his Theorem~13.1),
so we give an independent proof.

\newtheorem{theorem1}[theorem0]{Theorem}
\begin{theorem1}
Let $Q(t)$ be irreducible 		%
and define a linear recurrence by
(\ref{eq:GLR}), with at least one of $x_0, \ldots, x_{r-1}$ odd.
Then the sequence $(x_n)$ has period
\[	p_w \le 2^{w-2}\lambda	\]
for all $w \ge 2$ if $Q(-t)$ satisfies Condition~S,
\[	p_w \le 2^{w-2}\lambda	\]
for all $w \ge 3$ if $Q(t)$ satisfies Condition~S, and
\[	p_w = 2^{w-1}\lambda		\]
for all $w \ge 1$ iff neither $Q(t)$ nor $Q(-t)$ satisfies Condition~S.
\label{thm:T1}
\end{theorem1}

\leftline{\bf Proof}
\medskip

From Lemma~\ref{lemma:L_same_p}, $p_w = \rho_w$
is the order of $t \bmod (2^w, Q(t))$.
If $Q(-t)$ satisfies Condition~S then, from Theorem~\ref{thm:T0},
\[	t^\lambda = 1	\mymod (4, Q(t)).	\]
Using~(\ref{eq:trivia1}), it follows by induction on $w$ that
\[	t^{2^{w-2}\lambda} = 1	\mymod (2^w, Q(t))	\]
for all $w \ge 2$. This proves the first part of the
Theorem.	%
The second part is similar, so it only remains to prove the third part.

Suppose that $\rho_w = 2^{w-1}\lambda$ for all $w > 0$.
In particular, for $w = 3$ we have period $\rho_3 = 4\lambda$.
Thus
\[	t^{2\lambda} \ne 1 \mymod (8, Q(t))	\]
and, from~(\ref{eq:trivia3}),
\begin{equation}
	t^\lambda \ne \pm 1 \mymod  (4, Q(t)) . \label{eq:E4}
\end{equation}
From Theorem~\ref{thm:T0}, neither $Q(t)$
nor $Q(-t)$ can satisfy Condition~S, or we would obtain a
contradiction to~(\ref{eq:E4}).

Conversely, if neither $Q(t)$ or $Q(-t)$ satisfies Condition~S,
then we show by induction on $w$ that
\begin{equation}
	t^{2^{w-1}\lambda} = 1 + 2^wR_w	\mymod Q(t) ,	\label{eq:ind_w}
\end{equation}
where
\begin{equation}
	R_w \ne 0	\mymod (2, Q(t)) ,		\label{eq:R_w_odd}
\end{equation}
for all $w \ge 1$. Certainly
\[	t^\lambda = 1 \mymod (2, Q(t))	\]
but, from Theorem~\ref{thm:T0},
\[	t^\lambda \ne 1 \mymod (4, Q(t)) ,	\]
so (\ref{eq:ind_w}) and (\ref{eq:R_w_odd}) hold for $w = 1$.
Defining
\begin{equation}
	R_w = R_{w-1}(1 + 2^{w-2}R_{w-1})		\label{eq:R_rec}
\end{equation}
for $w \ge 2$, we see that (\ref{eq:ind_w}) holds for all $w \ge 1$.
It remains to prove (\ref{eq:R_w_odd}) for $w > 1$.
For $w = 2$, (\ref{eq:R_w_odd})
follows from Theorem~\ref{thm:T0} and~(\ref{eq:trivia3}), because
$t^\lambda \ne \pm 1 \bmod (4, Q(t))$
implies $t^{2\lambda} \ne 1 \bmod (8, Q(t))$.
For $w > 2$, (\ref{eq:R_w_odd}) follows by
induction from (\ref{eq:R_rec}), since $2^{w-2}$ is even.
It follows that
\[	\rho_w = 2^{w-1}\lambda \]
for all $w \ge 1$.						\QED

\section{Primitive Trinomials}
\label{sec:prim_tri}

In this section we consider a case of interest
because of its applications to pseudo-random number generation:
\[	Q(t) = q_0 + q_st^s + q_rt^r	\]
is a trinomial ($r > s > 0$).
Theorem~\ref{thm:T2} shows that the period is always maximal in
cases of practical interest.
The condition $r > 2$ is necessary,
as the example $Q(t) = 1 - t + t^2$ of Section~\ref{sec:cond_S} shows.

\newtheorem{theorem2}[theorem0]{Theorem}
\begin{theorem2}
Let $Q(t)$ be a primitive trinomial of degree $r > 2$.
Then the sequence $(x_n)$ defined by~(\ref{eq:GLR}) (with at least one
of $x_0, \ldots, x_{r-1}$ odd)
has period $p_w = 2^{w-1}(2^r - 1)$.
\label{thm:T2}
\end{theorem2}

\leftline{\bf Proof}
\medskip

From Theorem~\ref{thm:T1} it is sufficient to show that $Q(t)$
does not satisfy Condition~S.
(Since $Q(-t)$ is also a trinomial, the same argument shows that
$Q(-t)$ does not satisfy Condition~S.)

Suppose, by way of contradiction, that $Q(t)$ satisfies Condition~S.
We use the formulation of Condition~S given in Lemma~\ref{lemma:L_cond_S}.
Since $Q(t)$ is irreducible, $q_0 = q_s = q_r = 1 \bmod 2$.
If $s$ is even, say $s = 2m$, then
\[ \sum_{\scriptstyle j+k=2m \atop \scriptstyle 0 \le j < k \le r} q_jq_k
	= q_0q_s = 1	\mymod 2, \]
so $\epsilon_m \ne 0$, and (\ref{eq:epsdef}) implies that $q_m \ne 0$.
Since $0 < m < s < r$,
this contradicts the assumption that $Q(t)$ is a trinomial.
Hence, $s$ must be odd.

If $r$ is odd then
$r+s$ is even, and a similar argument shows that $q_{(r+s)/2} \ne 0$,
contradicting the assumption that $Q(t)$ is a trinomial.
Hence, $r$ must be even.

Taking $m = r/2$, we see that
$\epsilon_m \ne 0$, so $q_m \ne 0$.  This is only possible if
$m = s$, so
\[	Q(t) = t^{2s} + t^s + 1 \mymod 2 . \]
In this case $t^{3s} = 1 \bmod (2,Q(t))$.
Now $r = 2s > 2$, so $3s < 2^r - 1$, and $Q(t)$ can not be primitive.
This contradiction completes the proof.			\QED

A minor modification of the proof of Theorem~\ref{thm:T2} gives:

\newtheorem{theorem3}[theorem0]{Theorem}
\begin{theorem3}
Let $Q(t) = q_0 + q_st^s + q_rt^r$ be an irreducible trinomial
of degree $r \ne 2s$. Then
the sequence $(x_n)$ defined by~(\ref{eq:GLR}) (with at least one
of $x_0, \ldots, x_{r-1}$ odd)
has period $p_w = 2^{w-1}\lambda$.
\label{thm:T3}
\end{theorem3}

As mentioned above, it is easy to find primitive trinomials
of very high degree $r$ if $2^r - 1$ is a Mersenne prime.
Zierler~\cite{Zie69b} gives examples with
$r \le 9689$, and we found two examples with
higher degree: $t^{19937} + t^{9842} + 1$ and $t^{23209} + t^{9739} + 1$.
These and other examples with $r \le 44497$ were found
independently by Kurita and Matsumoto~\cite{Kur91}.
Such primitive trinomials provide the basis for fast random number
generators with extremely long periods
and good statistical properties~\cite{Bre92b}.

\section{Exceptional Polynomials}

We say that a polynomial $Q(t)$ of degree $r > 1$ is
{\em exceptional} if conditions 1--3 hold
and is a {\em candidate} if conditions 2--3 hold~--
\begin{enumerate}
\item $Q(t) \bmod 2$ is primitive.
\item $Q(t)$ has coefficients $q_j \in \{0, -1, +1 \}$, and $q_0 = q_r = 1$.
\item $Q(t)$ satisfies Condition~S.
\end{enumerate}

By Theorem~\ref{thm:T1}, if $Q(t)$ is exceptional then $Q(t)$
and $Q(-t)$ define simple linear recurrences (mod $2^w$)
which have less than the maximal period for $w > 2$.

Only the coefficients of $Q(t) \bmod 4$ are relevant to Condition~S.
If condition 2 is relaxed to allow coefficients equal to 2
then, by Lemma~\ref{lemma:L_cond_S},
there is one such $Q(t)$ corresponding to each primitive
polynomial in $Z_2[t]$.
With condition 2 as stated the number of these $Q(t)$ is considerably reduced.

It is interesting to consider strengthening condition 2 by asking
for certain patterns in the signs of the coefficients.
For example, we might ask for polynomials $Q(t)$ with all coefficients
$q_j \in \{0, 1\}$,
or for all coefficients of $\pm Q(-t)$ to be in $\{0, 1\}$.
There are candidates satisfying these conditions, but we have not
found any which are also exceptional, apart from the
trivial $Q(t) = 1 - t + t^2$.
It is possible for an exceptional polynomial to have
$(-1)^{j}q_j \ge 0$ for $0 \le j < r$.
The only example for $2 < r \le 44$ is	%
\[	Q(t) = 1 - t + t^2 - t^5 + t^6 + t^8 - t^9 +
		t^{10} + t^{12} - t^{13} + t^{16} + t^{18} + t^{21} . \]
Observe that
$Q(-t)$ defines a linear recurrence with nonnegative coefficients
\[ x_{n+21} = x_n + x_{n+1} + x_{n+2} + x_{n+5} + x_{n+6} + x_{n+8}
	+ x_{n+9} + x_{n+10} + x_{n+12} + x_{n+13} + x_{n+16} + x_{n+18} \]
which has period $p_2 = p_1 = 2^{21} - 1$ when considered mod~2 {\em or}
mod~4.

In Table~\ref{Tab:exceptions} we list the exceptional
polynomials $Q(t)$ of degree $r \le 14$.
If $Q(t)$ is exceptional then so is $\tilde Q(t)$.
Thus, we only list one of these in Table~\ref{Tab:exceptions}.

The number $\nu(r)$ of exceptional $Q(t)$
(counting only one of $Q(t), \tilde Q(t)$)
is given in Table~\ref{Tab:exception_kt}.
The term ``exceptional'' is
justified as $\nu(r)$ appears to
be a much more slowly growing function of $r$ than the number~\cite{Gol67}
\[	\lambda_2(r) = \varphi(2^r - 1)/r	\]
of primitive polynomials of degree $r$ in $Z_2[t]$
(where $\varphi$ is Euler's totient-function)
or the total number of polynomials of degree $r$ with coefficients
in $\{0, -1, +1 \}$.
Heuristic arguments suggest that
the number $\kappa(r)$ of candidates should grow like $(3/2)^r$ and that
$\nu(r)$ should grow like $(3/4)^r\lambda_2(r)$.
The arguments are as follows~--

\begin{quote}
There are $2^{r-1}$ polynomials ${\bar Q}(t)$
of degree $r$ with coefficients in
$\{0, 1\}$, satisfying ${\bar q}_0 = {\bar q}_r = 1$.
Randomly select such a ${\bar Q}(t)$, and compute
$\epsilon_0, \epsilon_1, \ldots, \epsilon_r$ from	%
\[ \sum_{\scriptstyle j+k=2m \atop \scriptstyle 0 \le j < k \le r}
	{\bar q}_j{\bar q}_k = \epsilon_m	\mymod 2 \]
Extend ${\bar Q}(t)$ to a polynomial $Q(t)$ with coefficients
$q_m \in \{-1, 0, 1, 2\}$ such that

${\bar q}_m = q_m \bmod 2$
and~(\ref{eq:epsdef}) is satisfied for $0 \le m \le r$.
The (unique) mapping is
given by

$q_m = {\bar q}_m + 2\epsilon_m \bmod 4$.
It is easy to see that $q_0 = q_r = 1$. If we {\em assume} that
each $q_m$ for $1 \le m < r$ has independent probability 1/4 of
assuming the ``forbidden'' value 2, then the probability that
$Q(t)$ is a candidate is $(3/4)^{r-1}$.
Thus,
\[	\kappa(r) \simeq (3/2)^{r-1}. \]
The argument is not strictly correct. For example,
it gives a positive probability that $q_1 = 0$, $q_2 = 1$, but this
never occurs for $r > 2$.  However, the argument does appear to predict
the correct order of magnitude of $\kappa(r)$.

The probability that a randomly chosen ${\bar Q}(t)$
with ${\bar q}_0 = {\bar q}_r = 1$ is primitive is
just $\lambda_2(r)/2^{r-1}$.  If there is the same probability that
a randomly chosen candidate is primitive, then the number of
primitive candidates should be $(3/4)^{r-1}\lambda_2(r)$,
and $\nu(r)$ should be half this number.
\end{quote}

In Table~\ref{Tab:exception_kt} we give
\[	\bar \nu(r) = {\nu(r) \over (3/4)^r\lambda_2(r)}\;\;; \]
the numerical evidence suggests that $\bar \nu(r)$ converges to a
positive constant $\bar \nu(\infty)$ as $r \to \infty$.
However, $\bar \nu(\infty)$ is less than the value 2/3
predicted by the heuristic argument.
Our best estimate (obtained from a separate computation which
gives faster convergence) is
\[	\bar \nu(\infty) = 0.45882 \pm 0.00002	\]

The computation of Table~\ref{Tab:exception_kt} took 166 hours
on a VaxStation 3100. 		%
We outline the method used.
It is easy to check if a candidate polynomial is exceptional~\cite{Knu81}.
A straightforward method of enumerating all
candidate polynomials of degree $r$ is to associate a
polynomial $Q(t)$ such that $q_0 = q_r = 1$
with an $(r-1)$-bit binary number $N = b_1 \cdots b_{r-1}$,
where $b_j = q_j \bmod 2$.  For each such $N$, compute
$\epsilon_0, \ldots, \epsilon_r$ from~(\ref{eq:epssum}).
Now~(\ref{eq:epsdef}) defines $q_0, \ldots, q_r \bmod 4$.
If there is an index $m$ such that $\epsilon_m = 1 \bmod 2$
but $q_m = 0 \bmod 2$, then~(\ref{eq:epsdef}) shows that
$q_m = 2 \bmod 4$, contradicting condition 2.
The straightforward enumeration has complexity $\Omega(2^r)$,
but this can be reduced by two devices~--
\begin{enumerate}
\item 	If~(\ref{eq:epsdef}) shows that $q_m = 2 \bmod 4$ for
        some $m < r/2$, we may use the fact that $\epsilon_m$
	in~(\ref{eq:epssum}) depends only on $q_0,\ldots,q_{2m}$
	to skip over a block of $2^{r-2m-1}$ numbers $N$.
	By an argument similar to the heuristic argument for
	the order of magnitude of $\nu(r)$, with support from
	empirical evidence for $r \le 40$,
	we conjecture that this device reduces the complexity
	of the enumeration to
\[	O\left({r^2}2^r\left({{3}\over{4}}\right)^{r/2}\right)
	= O({r^2}3^{r/2}) .	\]
\item	Fix $s$, $0 \le s < r$. Since $\epsilon_{r-m}$
	in~(\ref{eq:epssum}) depends only on $q_{r-2m},\ldots,q_r$,
	we can tabulate those low-order bits
	$b_{r-s} \cdots b_{r-1}$ which do not necessarily
	lead to condition 2 being violated for some $q_{r-m}$,
	$2m \le s$. In the enumeration we need only consider $N$ with
 	low-order bits in the table.
	We conjecture that this reduces the complexity
	of the enumeration to
\[	O\left({r^2}2^r\left({{3}\over{4}}\right)^{s/2}\right)
	= O({r^2}2^{r-s}3^{s/2}) 	\]
	provided care is taken to generate the table efficiently.
\end{enumerate}

The two devices can be combined, but they are not independent.
The complexity of the combination is conjectured to be
\[	O\left({r^2}2^r\left({{3}\over{4}}\right)^{(6r+5s)/12}\right)
	= O\left({r^2}3^{r/2}\left({{3}\over{4}}\right)^{5s/12}\right) ,\]
where the exponent $5s/12$ (instead of $s/2$) reflects the lack of
independence.
In the computation of Table~\ref{Tab:exception_kt}
we used $s \le 22$ because of memory constraints.
The table size is $O(s3^{s/2})$ bits if the table is stored
as a list to take advantage of sparsity.

\medskip
\leftline{\bf Acknowledgements}
\medskip

We thank a referee for pointing out an error in the formulation
of Lemma~\ref{lemma:L_same_p} given in~\cite{Bre92a}, and for providing
references to the classical literature.
The ANU Supercomputer Facility provided time on a Fujitsu VP~2200/10 for
the discovery of the primitive
trinomials mentioned at the end of Section~\ref{sec:prim_tri}.

\begin{table}[p]
\renewcommand{\baselinestretch}{1.0}	%
\centerline{
\begin{tabular}{|c|c|} \hline
$r$     & $Q(t)$ \\ \hline
2 	& $1-t+t^2$\\ \hline
5 	& $1-t-t^2+t^4+t^5$\\ \hline
 	& $1-t+t^2+t^3-t^4-t^6+t^9$\\				%
9  	& $1-t+t^2-t^3-t^4+t^8+t^9$\\
  	& $1-t+t^2-t^3-t^4-t^5+t^6+t^8+t^9$\\ \hline
10 	& $1-t+t^2+t^3+t^4+t^6-t^7+t^9+t^{10}$\\ \hline
11 	& $1-t+t^2-t^3-t^4+t^5+t^6-t^8+t^{11}$\\ \hline
12 	& $1-t+t^2-t^3-t^4-t^8+t^9+t^{11}+t^{12}$\\ \hline
 	& $1-t+t^2-t^3+t^4-t^5-t^6+t^{12}+t^{13}$\\		%
   	& $1-t+t^2-t^3+t^4-t^5-t^6-t^7+t^8+t^{12}+t^{13}$\\
13   	& $1-t-t^2-t^4-t^6+t^7-t^8+t^9+t^{10}+t^{12}+t^{13}$\\
   	& $1-t+t^2+t^3+t^4+t^5+t^7+t^9-t^{11}-t^{12}+t^{13}$\\
   	& $1-t+t^2+t^3+t^4+t^5-t^8-t^9-t^{11}-t^{12}+t^{13}$\\ \hline
 	& $1-t+t^2+t^3-t^4-t^6-t^7+t^8+t^9-t^{11}+t^{14}$\\	%
   	& $1+t+t^3-t^4-t^5+t^6+t^7+t^8+t^9-t^{11}+t^{14}$\\
14   	& $1-t-t^2+t^3-t^5+t^6+t^7-t^8-t^9+t^{13}+t^{14}$\\
   	& $1-t-t^2-t^3-t^5+t^7+t^9+t^{10}-t^{11}+t^{13}+t^{14}$\\
   	& $1-t-t^2+t^4-t^6+t^8+t^9+t^{10}+t^{11}+t^{13}+t^{14}$\\ \hline
\end{tabular}
}
\caption{Exceptional Polynomials of degree $r \le 14$}
\label{Tab:exceptions}	%

\end{table}

\begin{table}[p]
\renewcommand{\baselinestretch}{1.0}	%
\centerline{
\begin{tabular}{|c|c|c||c|c|c|} \hline
$r$     & $\nu(r)$ & $\bar\nu(r)$ & $r$	& $\nu(r)$ & $\bar\nu(r)$\\ \hline
1	& 0	& 0	& 21	& 79	& 0.3923 \\
2	& 1	& 1.78	& 22	& 94	& 0.4390 \\
3	& 0	& 0	& 23	& 231	& 0.4837 \\
4	& 0	& 0	& 24 	& 129	& 0.4650 \\
5	& 1	& 0.70	& 25	& 428	& 0.4388 \\
6	& 0	& 0	& 26	& 448	& 0.4615 \\
7	& 0	& 0	& 27	& 883	& 0.4964 \\
8	& 0	& 0	& 28	& 635	& 0.4218 \\
9	& 3	& 0.83	& 29	& 1933	& 0.4410 \\
10	& 1	& 0.30	& 30	& 1470	& 0.4619 \\
11	& 1	& 0.13	& 31	& 4380	& 0.4721 \\
12	& 1	& 0.22	& 32	& 3125	& 0.4636 \\
13	& 5	& 0.33	& 33	& 7232	& 0.4549 \\
14	& 5	& 0.37	& 34	& 8862	& 0.4656 \\
15	& 15	& 0.62	& 35	& 18870 & 0.4792 \\
16	& 12	& 0.58	& 36	& 10516 & 0.4560 \\
17	& 26	& 0.45	& 37	& 40082	& 0.4547 \\
18	& 18	& 0.41	& 38	& 39858	& 0.4623 \\
19	& 62	& 0.53	& 39	& 75370	& 0.4712 \\
20	& 34	& 0.45	& 40	& 54758	& 0.4598 \\
\hline
\end{tabular}
}
\caption{Number of Exceptional Polynomials}
\label{Tab:exception_kt}     %
\end{table}

\end{document}